\documentclass[11pt]{amsart}
\usepackage{mathrsfs}
\newtheorem{theorem}{Theorem}[section]
\newtheorem{proposition}[theorem]{Proposition}
\newtheorem{lemma}[theorem]{Lemma}
\newtheorem{corollary}[theorem]{Corollary}
\theoremstyle{definition}
\newtheorem{definition}[theorem]{Definition}

\begin{document}

\title[Rotational symmetry of ancient solutions to the Ricci flow]{Rotational symmetry of ancient solutions to the Ricci flow in dimension $3$ -- The compact case}
\author{Simon Brendle}
\address{Department of Mathematics, Columbia University, 2990 Broadway, New York, NY 10027, USA}
\begin{abstract}
In \cite{Brendle2}, we proved that every noncompact ancient $\kappa$-solution to the Ricci flow in dimension $3$ is either locally isometric to a family of shrinking cylinders, or isometric to the Bryant soliton. In the same paper, we announced that the same method implies that compact ancient $\kappa$-solutions are rotationally symmetric. In this note we provide the details of this argument.
\end{abstract}
\thanks{The author was supported by the National Science Foundation under grant DMS-1806190 and by the Simons Foundation.}
\maketitle

\section{Introduction}

In a recent paper \cite{Brendle2}, we proved:

\begin{theorem}
\label{classification.noncompact.case}
Assume that $(M,g(t))$ is a three-dimensional ancient $\kappa$-solution which is noncompact. Then $(M,g(t))$ is isometric to either a family of shrinking cylinders (or a quotient thereof), or to the Bryant soliton. 
\end{theorem}

Theorem \ref{classification.noncompact.case} confirms a conjecture of Perelman \cite{Perelman1}.

The proof of Theorem \ref{classification.noncompact.case} relies on a number of crucial ingredients: first, Perelman's Canonical Neighborhood Theorem for ancient $\kappa$-solutions; second, an estimate due to Anderson-Chow \cite{Anderson-Chow} for the parabolic Lichnerowicz equation; third, the classification of $\kappa$-noncollapsed steady gradient Ricci solitons (cf. \cite{Brendle1}); fourth, the classification of ancient $\kappa$-solutions with rotational symmetry (cf. \cite{Brendle2}, Part 1); and fifth, the Neck Improvement Theorem, which asserts that a neck becomes more and more rotationally symmetric under the evolution (cf. \cite{Brendle2}, Part 2).

As we announced in our November 6, 2018 preprint, the arguments in that paper also imply that compact ancient $\kappa$-solutions are rotationally symmetric:

\begin{theorem}
\label{main.thm}
Let $(M,g(t))$ be a three-dimensional ancient $\kappa$-solution which is compact and simply connected. Then $(M,g(t))$ is rotationally symmetric.
\end{theorem}

In this short note, we will explain how Theorem \ref{main.thm} follows from the results and techniques developed in \cite{Brendle2}. 

After the proof of Theorem \ref{main.thm} was written, we have learned of an alternative approach by Bamler-Kleiner \cite{Bamler-Kleiner}. 

\section{Structure of compact ancient $\kappa$-solutions}

In this section, we recall some basic facts about the structure of compact ancient $\kappa$-solutions. Throughout this section, we assume that $(M,g(t))$ is a three-dimensional ancient $\kappa$-solution which is compact and simply connected. Moreover, we assume that $(M,g(t))$ is not a family of shrinking round spheres. Roughly speaking, if $-t$ is large, then the manifold $(M,g(t))$ will have the structure of a tube with two caps attached on each side, and each of these caps look like a piece of the Bryant soliton. In the following, we make this precise.

\begin{proposition}
\label{asymptotic.soliton}
The asymptotic shrinking soliton associated with $(M,g(t))$ is a cylinder.
\end{proposition}

\textbf{Proof.} 
By Perelman's classification of shrinking gradient Ricci solitons in dimension $3$ (cf. \cite{Perelman2}), the asymptotic shrinking soliton associated with $(M,g(t))$ is either a sphere or a cylinder. If the asymptotic shrinking soliton associated with $(M,g(t))$ is a sphere, then, by Hamilton's curvature pinching estimates, the solution $(M,g(t))$ has constant sectional curvature for each $t$, contrary to our assumption. \\

\begin{proposition}
\label{possible.limit.flows}
Let $(x_k,t_k)$ be an arbitrary sequence of points in space-time satisfying $\lim_{k \to \infty} t_k = -\infty$. Let us perform a parabolic rescaling around the point $(x_k,t_k)$ by the factor $R(x_k,t_k)$. After passing to a subsequence, the rescaled flows converge to a limit which is either a family of shrinking cylinders or the Bryant soliton. 
\end{proposition}

\textbf{Proof.} 
By Perelman's work \cite{Perelman1}, the rescaled manifolds converge to an ancient $\kappa$-solution. If the limiting ancient solution is noncompact, then, by \cite{Brendle2}, it must be either a family of shrinking cylinders or the Bryant soliton, and we are done. Hence, it remains to consider the case when the limiting ancient solution is compact. In this case, we have $\limsup_{k \to \infty} R(x_k,t_k) \, \text{\rm diam}_{g(t_k)}(M)^2 < \infty$ and furthermore $\limsup_{k \to \infty} R(x_k,t_k)^{-1} \, \sup_M R(x,t_k) < \infty$. Putting these facts together gives $\limsup_{k \to \infty} \text{\rm diam}_{g(t_k)}(M)^2 \, \sup_M R(x,t_k) < \infty$. This implies that $(M,g(t_k))$ cannot contain arbitrarily long necks. On the other hand, since the asymptotic shrinking soliton is a cylinder by Proposition \ref{asymptotic.soliton}, we know that $(M,g(t_k))$ must contain arbitrarily long necks if $k$ is sufficiently large. This is a contradiction. \\

In the next step, we fix a small number $\varepsilon_1>0$. For later purposes, it is important that we choose $\varepsilon_1$ small enough so that the Neck Improvement Theorem in \cite{Brendle2} holds. Let us fix a small number $\theta>0$ with the following property: if $(x,t)$ is a point in space-time satisfying $\lambda_1(x,t) \leq \theta R(x,t)$, then the point $(x,t)$ lies at the center of an evolving $\varepsilon_1$-neck. \\

\begin{proposition}
\label{caps}
If $-t$ is sufficiently large, then we can find two disjoint domains $D_1$ and $D_2$ (depending on $t$) with the following properties: 
\begin{itemize}
\item For each point $x \in M \setminus (D_1 \cup D_2)$, we have $\lambda_1(x,t) \leq \theta R(x,t)$. In particular, the point $(x,t)$ lies at the center of an evolving $\varepsilon_1$-neck.
\item $D_1$ and $D_2$ are diffeomorphic to $B^3$.
\item $\partial D_1$ and $\partial D_2$ are leafs of Hamilton's CMC foliation.
\item There exists a point $q_1 \in \partial D_1$ where $\lambda_1(q_1,t) \geq \frac{2}{3} \, \theta R(q_1,t)$.
\item There exists a point $q_2 \in \partial D_2$ where $\lambda_1(q_2,t) \geq \frac{2}{3} \, \theta R(q_2,t)$.
\end{itemize}
\end{proposition}

\textbf{Proof.} 
By Proposition \ref{asymptotic.soliton}, the asymptotic soliton is a cylinder. Hence, if $-t$ is large, there is a point $q \in M$ (depending on $t$) where $\lambda_1(q,t) \leq \frac{1}{2} \, \theta R(q,t)$. In particular, $(q,t)$ lies at the center of an $\varepsilon_1$-neck. We follow this neck to each side until we encounter a point with $\lambda_1 \geq \frac{3}{4} \, \theta R$. From this, the assertion follows. \\

\begin{proposition}
\label{caps.look.like.bryant.soliton}
Consider a sequence of times $t_k \to -\infty$. Let $D_{1,k}$ and $D_{2,k}$ denote the domains defined in the previous proposition, and Let $q_{1,k}$ and $q_{2,k}$ be defined as in the previous proposition. If we rescale around the point $(q_{1,k},t_k)$, then the rescaled flows converge in the Cheeger-Gromov sense to the Bryant soliton. Similarly, if we rescale around the point $(q_{2,k},t_k)$, then the rescaled flows converge in the Cheeger-Gromov sense to the Bryant soliton. 
\end{proposition}

\textbf{Proof.} 
We know that the rescaling limit is either a family of shrinking cylinders or the Bryant soliton. Since $\lambda_1(q_{1,k},t_k) \geq \frac{1}{2} \, \theta R(q_{1,k},t_k)$ and $\lambda_1(q_{2,k},t_k) \geq \frac{1}{2} \, \theta R(q_{2,k},t_k)$, the limit cannot be a cylinder, so it must be the Bryant soliton. \\

\begin{corollary}
If $-t$ is sufficiently large, then $\lambda_1(x,t) \geq \frac{1}{2} \, \theta R(x,t)$ for each $x \in D_1 \cup D_2$.
\end{corollary}

\textbf{Proof.} 
If $-t$ is sufficiently large, the domain $D_1$ will be arbitrarily close to a subset of the Bryant soliton after rescaling. Since $\lambda_1 \geq \frac{2}{3} \, \theta R$ at some point on $\partial D_1$, it follows that $\lambda_1 \geq \frac{1}{2} \, \theta R$ in $D_1$. An analogous statement gives $\lambda_2 \geq \frac{1}{2} \, \theta R$ in $D_2$. \\

\begin{corollary}
\label{local.maxima.of.scalar.curvature}
If $-t$ is sufficiently large, then there exists a unique point $p_1 \in D_1$ (depending on $t$) such that $R(p_1,t) = \sup_{x \in D_1} R(x,t)$. Similarly, there exists a unique point $p_2 \in D_2$ (depending on $t$) such that $R(p_2,t) = \sup_{x \in D_1} R(x,t)$.
\end{corollary}

\begin{proposition}
\label{distance.between.caps}
Suppose we are given a large constant $Q$. Suppose that $-t$ is sufficiently large (depending in $Q$), and let $p_1,p_2$ be defined as in Corollary \ref{local.maxima.of.scalar.curvature}. Then the balls $B_{g(t)}(p_1,Q \, R(p_1,t)^{-\frac{1}{2}})$ and $B_{g(t)}(p_2,Q \, R(p_2,t)^{-\frac{1}{2}})$ are disjoint. 
\end{proposition}

\textbf{Proof.} 
There exists a small constant $\beta>0$ such that $\lambda_1(x,t) \geq \beta R(x,t)$ for all $x \in B_{g(t)}(p_1,Q \, R(p_1,t)^{-\frac{1}{2}}) \cup B_{g(t)}(p_2,Q \, R(p_2,t)^{-\frac{1}{2}})$. On the other hand, Proposition \ref{asymptotic.soliton} implies that, for $-t$ sufficiently large, we can find a neck $N \subset M \setminus (D_1 \cup D_2)$ with the property that $\lambda_1(x,t) < \beta R(x,t)$ for all $x \in N$. This gives $B_{g(t)}(p_1,Q \, R(p_1,t)^{-\frac{1}{2}}) \cap N = B_{g(t)}(p_2,Q \, R(p_2,t)^{-\frac{1}{2}}) \cap N = \emptyset$. The balls $B_{g(t)}(p_1,Q \, R(p_1,t)^{-\frac{1}{2}})$ and $B_{g(t)}(p_2,Q \, R(p_2,t)^{-\frac{1}{2}})$ lie in different connected components of $M \setminus N$, and hence must be disjoint. \\

\begin{proposition}
\label{points.away.from.tip.are.necklike}
Let $\varepsilon>0$ be given. Then there exists a large constant $Q$ (depending on $\varepsilon$) with the following property. Suppose that $-t$ is sufficiently large, and let $p_1,p_2$ be defined as in Corollary \ref{local.maxima.of.scalar.curvature}. Moreover, suppose that $x$ is a point satisfying $d_{g(t)}(p_1,x) \geq Q \, R(p_1,t)^{-\frac{1}{2}}$ and $d_{g(t)}(p_2,x) \geq Q \, R(p_2,t)^{-\frac{1}{2}}$. Then $(x,t)$ lies on an evolving $\varepsilon$-neck. 
\end{proposition}

\textbf{Proof.} 
Suppose this is false. Then there exists a sequence of times $t_k \to -\infty$ and a sequence of points $(x_k,t_k)$ in space-time such that $d_{g(t_k)}(p_{k,1},x_k)^2 \, R(p_{1,k},t_k) \to \infty$, $d_{g(t_k)}(p_{k,1},x_k)^2 \, R(p_{2,k},t_k) \to \infty$, and $(x_k,t_k)$ does not lie on an evolving $\varepsilon$-neck. Clearly, $x_k \in M \setminus (D_{1,k} \cup D_{2,k})$ for $k$ large. In particular, $(x_k,t_k)$ lies on an $\varepsilon_1$-neck. Let $\Sigma_k$ denote the leaf of the CMC foliation passing through $(x_k,t_k)$. We now rescale the flow around $(x_k,t_k)$ by the factor $R(x_k,t_k)$, and pass to the limit as $k \to \infty$. Since $(x_k,t_k)$ does not lie on an $\varepsilon$-neck, the limit cannot be a cylinder. Consequently, the limit must be the Bryant soliton. Therefore, we can find a domain $\Omega_k \subset M$ such that $\partial \Omega_k = \Sigma_k$, $\limsup_{k \to \infty} \text{\rm diam}_{g(t_k)}(\Omega_k)^2 \, \sup_{x \in \Omega_k} R(x,t_k) < \infty$. If $k$ is sufficiently large, we either have $p_{1,k} \in \Omega_k$ or $p_{2,k} \in \Omega_k$. If $p_{1,k} \in \Omega_k$, then $\liminf_{k \to \infty} d_{g(t_k)}(p_{1,k},x_k)^2 \, R(p_{1,k},t_k) < \infty$, contrary to our assumption. Similarly, if $p_{1,k} \in \Omega_k$, then $\liminf_{k \to \infty} d_{g(t_k)}(p_{2,k},x_k)^2 \, R(p_{2,k},t_k) < \infty$, which again contradicts our assumption. \\

\section{Proof of rotational symmetry}

In this section, we give the proof of rotational symmetry. Throughout this section, we assume that $(M,g(t))$ is a three-dimensional ancient $\kappa$-solution which is compact and simply connected. Moreover, we assume that $(M,g(t))$ is not a family of shrinking round spheres. We claim that $(M,g(t))$ is rotationally symmetric. The proof is by contradiction, and we will assume throughout this section that $(M,g(t))$ is not rotationally symmetric. 

We begin with a definition, which is adapted from \cite{Brendle2}:

\begin{definition}
\label{symmetry.of.cap}
We say that $(M,g(\bar{t}))$ is $\varepsilon$-symmetric if there exist a compact domain $D \subset M$ and time-independent vector fields $U^{(1)},U^{(2)},U^{(3)}$ on $D$ with the following properties: 
\begin{itemize} 
\item The domain $D$ is a disjoint union of two domains $D_1$ and $D_2$, each of which is diffeomorphic to $B^3$. 
\item $\partial D_1$ and $\partial D_2$ are leaves of Hamilton's CMC foliation of $(M,g(\bar{t}))$. We put $\rho_1^{-2} := \sup_{x \in D_1} R(x,\bar{t})$ and $\rho_2^{-2} := \sup_{x \in D_2} R(x,\bar{t})$.
\item $\lambda_1(x,\bar{t}) \leq \theta R(x,\bar{t})$ for all points $x \in M \setminus D$.
\item $\lambda_1(x,\bar{t}) \geq \frac{1}{2} \, \theta R(x,\bar{t})$ for all points $x \in D$.
\item For each $x \in M \setminus D$, the point $(x,\bar{t})$ is $\varepsilon$-symmetric in the sense of \cite{Brendle2}.
\item $\sup_{D_1 \times [\bar{t}-\rho_1^2,\bar{t}]} \sum_{l=0}^2 \sum_{a=1}^3 \rho_1^{2l} \, |D^l(\mathscr{L}_{U^{(a)}}(g(t)))|^2 \leq \varepsilon^2$.
\item $\sup_{D_2 \times [\bar{t}-\rho_2^2,\bar{t}]} \sum_{l=0}^2 \sum_{a=1}^3 \rho_2^{2l} \, |D^l(\mathscr{L}_{U^{(a)}}(g(t)))|^2 \leq \varepsilon^2$.
\item If $\Sigma \subset D_1$ is a leaf of the CMC foliation of $(M,g(\bar{t}))$ which has distance at most $\text{\rm area}_{g(\bar{t})}(\partial D_1)^{\frac{1}{2}}$ from $\partial D_1$, then $\sup_\Sigma \sum_{a=1}^3 \rho_1^{-2} \, |\langle U^{(a)},\nu \rangle|^2 \leq \varepsilon^2$, where $\nu$ denotes the unit normal vector to $\Sigma$ in $(M,g(\bar{t}))$.
\item If $\Sigma \subset D_2$ is a leaf of the CMC foliation of $(M,g(\bar{t}))$ which has distance at most $\text{\rm area}_{g(\bar{t})}(\partial D_2)^{\frac{1}{2}}$ from $\partial D_2$, then $\sup_\Sigma \sum_{a=1}^3 \rho_2^{-2} \, |\langle U^{(a)},\nu \rangle|^2 \leq \varepsilon^2$, where $\nu$ denotes the unit normal vector to $\Sigma$ in $(M,g(\bar{t}))$.
\item If $\Sigma \subset D_1$ is a leaf of the CMC foliation of $(M,g(\bar{t}))$ which has distance at most $\text{\rm area}_{g(\bar{t})}(\partial D_1)^{\frac{1}{2}}$ from $\partial D_1$, then 
\[\sum_{a,b=1}^3 \bigg | \delta_{ab} - \text{\rm area}_{g(\bar{t})}(\Sigma)^{-2} \int_\Sigma \langle U^{(a)},U^{(b)} \rangle \, d\mu_{g(\bar{t})} \bigg |^2 \leq \varepsilon^2.\]
\item If $\Sigma \subset D_2$ is a leaf of the CMC foliation of $(M,g(\bar{t}))$ which has distance at most $\text{\rm area}_{g(\bar{t})}(\partial D_2)^{\frac{1}{2}}$ from $\partial D_2$, then 
\[\sum_{a,b=1}^3 \bigg | \delta_{ab} - \text{\rm area}_{g(\bar{t})}(\Sigma)^{-2} \int_\Sigma \langle U^{(a)},U^{(b)} \rangle \, d\mu_{g(\bar{t})} \bigg |^2 \leq \varepsilon^2.\]
\end{itemize}
\end{definition}

\begin{proposition}
\label{solution.becomes.more.and.more.symmetric.as.t.approaches.minus.infinity}
Let $\varepsilon>0$ be given. If $-t$ is sufficiently large (depending on $\varepsilon$), then $(M,g(t))$ is $\varepsilon$-symmetric.
\end{proposition}

\textbf{Proof.} 
This follows from Proposition \ref{caps.look.like.bryant.soliton} together with Proposition \ref{points.away.from.tip.are.necklike}. \\

We next consider an arbitrary sequence $\varepsilon_k \to 0$. For $k$ large, we define 
\[t_k = \inf \{t \in (-\infty,0]: \text{\rm $(M,g(t))$ is not $\varepsilon_k$-symmetric}\}.\]  
Clearly, $(M,g(t))$ is $\varepsilon_k$-symmetric for each $t \in (-\infty,t_k)$. If $\limsup_{k \to \infty} t_k > -\infty$, then we conclude that $(M,g(t))$ is rotationally symmetric for $-t$ sufficiently large, and this contradicts our assumption. Therefore, $\limsup_{k \to \infty} t_k = -\infty$. 

\begin{lemma}
\label{varepsilon_k.symmetry.at.earlier.times}
If $t \in (-\infty,t_k)$, then $(M,g(t))$ is $\varepsilon_k$-symmetric. In particular, if $(x,t) \in M \times (-\infty,t_k)$ is a point in spacetime satisfying $\lambda_1(x,t) < \frac{1}{2} \theta R(x,t)$, then the point $(x,t)$ is $\varepsilon_k$-symmetric.
\end{lemma}

\textbf{Proof.} 
The first statement follows directly from the definition of $t_k$. The second statement follows from Definition \ref{symmetry.of.cap}. \\

For each $t$, we denote by $p_{1,t}$ and $p_{2,t}$ the local maxima of scalar curvature constructed in Corollary \ref{local.maxima.of.scalar.curvature}. Let $L$ denote the constant in the Neck Improvement Theorem in \cite{Brendle2}. By Proposition \ref{points.away.from.tip.are.necklike}, we can find a time $T \in (-\infty,0]$ and a large constant $\Lambda$ with the following property: if $\bar{t} \leq T$ and $d_{g(\bar{t})}(p_{1,\bar{t}},x) \geq \Lambda \, R(p_{1,\bar{t}},\bar{t})^{-\frac{1}{2}}$ and $d_{g(\bar{t})}(p_{2,\bar{t}},x) \geq \Lambda \, R(p_{2,\bar{t}},\bar{t})^{-\frac{1}{2}}$, then $\lambda_1(x,t) < \frac{1}{2} \theta R(x,t)$ for all points $(x,t) \in B_{g(\bar{t})}(\bar{x},L \, R(\bar{x},\bar{t})^{-\frac{1}{2}}) \times [\bar{t} - L \, R(\bar{x},\bar{t})^{-1},\bar{t}]$.

\begin{lemma}
\label{improved.symmetry.for.faraway.points}
Suppose that $\bar{t} \in (-\infty,t_k]$ and $d_{g(\bar{t})}(p_{1,\bar{t}},x) \geq \Lambda \, R(p_{1,\bar{t}},\bar{t})^{-\frac{1}{2}}$ and $d_{g(\bar{t})}(p_{2,\bar{t}},x) \geq \Lambda \, R(p_{2,\bar{t}},\bar{t})^{-\frac{1}{2}}$. Then $(\bar{x},\bar{t})$ is $\frac{\varepsilon_k}{2}$-symmetric.
\end{lemma}

\textbf{Proof.} 
By our choice of $\Lambda$, every point in the parabolic neighborhood $B_{g(\bar{t})}(\bar{x},L \, R(\bar{x},\bar{t})^{-\frac{1}{2}}) \times [\bar{t} - L \, R(\bar{x},\bar{t})^{-1},\bar{t}]$ satisfies $\lambda_1(x,t) < \frac{1}{2} \theta R(x,t)$. By our choice of $\theta$, every point in the parabolic neighborhood $B_{g(\bar{t})}(\bar{x},L \, R(\bar{x},\bar{t})^{-\frac{1}{2}}) \times [\bar{t} - L \, R(\bar{x},\bar{t})^{-1},\bar{t}]$ lies at the center of an evolving $\varepsilon_1$-neck. Moreover, by Lemma \ref{varepsilon_k.symmetry.at.earlier.times}, every point in the parabolic neighborhood $B_{g(\bar{t})}(\bar{x},L \, R(\bar{x},\bar{t})^{-\frac{1}{2}}) \times [\bar{t} - L \, R(\bar{x},\bar{t})^{-1},\bar{t})$ is $\varepsilon$-symmetric. Hence, the Neck Improvement Theorem in \cite{Brendle2} implies that $(\bar{x},\bar{t})$ is $\frac{\varepsilon_k}{2}$-symmetric. \\

\begin{proposition}
\label{cap.is.close.to.bryant}
If we rescale the solution around $(p_{1,t_k},t_k)$ by the factor $\rho_{1,k}^{-2} := R(p_{1,t_k},t_k)$, then the rescaled flows converge to the Bryant soliton. Similarly, if we rescale the solution around $(p_{1,t_k},t_k)$ by the factor $\rho_{2,k}^{-2} := R(p_{2,t_k},t_k)$, then the rescaled flows converge to the Bryant soliton. 
\end{proposition}

\textbf{Proof.} 
Since we know that $t_k \to -\infty$, this follows from Proposition \ref{caps.look.like.bryant.soliton}. \\

\begin{corollary} 
\label{almost.Bryant}
There exists a sequence $\delta_k \to 0$ such that the following statements hold when $k$ is sufficiently large: 
\begin{itemize}
\item The points $p_{1,t}$ depends smoothly on $t$, and $|\frac{d}{dt} p_{1,t}|_{g(t)} \leq \delta_k \rho_{1,k}^{-1}$ for $t \in [t_k-\delta_k^{-1} \rho_{1,k}^2,t_k]$.
\item The points $p_{2,t}$ depends smoothly on $t$, and $|\frac{d}{dt} p_{2,t}|_{g(t)} \leq \delta_k \rho_{2,k}^{-1}$ for $t \in [t_k-\delta_k^{-1} \rho_{2,k}^2,t_k]$.
\item The scalar curvature satisfies $\frac{1}{2K} \, (\rho_{1,k}^{-1} \, d_{g(t)}(p_{1,t},x)+1)^{-1} \leq \rho_{1,k}^2 \, R(x,t) \leq 2K \, (\rho_{1,k}^{-1} \, d_{g({1,t})}(p_{1,t},x)+1)^{-1}$ for all points $(x,t) \in B_{g(t_k)}(p_{1,t_k},\delta_k^{-1} \rho_{1,k}) \times [t_k-\delta_k^{-1} \rho_{1,k}^2,t_k]$.
\item The scalar curvature satisfies $\frac{1}{2K} \, (\rho_{2,k}^{-1} \, d_{g(t)}(p_{2,t},x)+1)^{-1} \leq \rho_{2,k}^2 \, R(x,t) \leq 2K \, (\rho_{2,k}^{-1} \, d_{g({2,t})}(p_{2,t},x)+1)^{-1}$ for all points $(x,t) \in B_{g(t_k)}(p_{2,t_k},\delta_k^{-1} \rho_{2,k}) \times [t_k-\delta_k^{-1} \rho_{2,k}^2,t_k]$.
\item There exists a nonnegative function $f: B_{g(t_k)}(p_{1,t_k},\delta_k^{-1} \rho_{1,k}) \times [t_k-\delta_k^{-1} \rho_{1,k}^2,t_k] \to \mathbb{R}$ such that $|\text{\rm Ric}-D^2 f| \leq \delta_k \rho_{1,k}^{-2}$, $|\Delta f + |\nabla f|^2 - \rho_{1,k}^{-2}| \leq \delta_k \rho_{1,k}^{-2}$, and $|\frac{\partial}{\partial t} f + |\nabla f|^2| \leq \delta_k \rho_{1,k}^{-2}$. Moreover, the function $f$ satisfies $\frac{1}{2K} \, (\rho_{1,k}^{-1} \, d_{g(t)}(p_{1,t},x)+1) \leq f(x,t)+1 \leq 2K \, (\rho_{1,k}^{-1} \, d_{g(t)}(p_{1,t},x)+1)$ for all points $(x,t) \in B_{g(t_k)}(p_{1,t_k},\delta_k^{-1} \rho_{1,k}) \times [t_k-\delta_k^{-1} \rho_{1,k}^2,t_k]$.
\item There exists a nonnegative function $f: B_{g(t_k)}(p_{2,t_k},\delta_k^{-1} \rho_{2,k}) \times [t_k-\delta_k^{-1} \rho_{2,k}^2,t_k] \to \mathbb{R}$ such that $|\text{\rm Ric}-D^2 f| \leq \delta_k \rho_{2,k}^{-2}$, $|\Delta f + |\nabla f|^2 - \rho_{2,k}^{-2}| \leq \delta_k \rho_{2,k}^{-2}$, and $|\frac{\partial}{\partial t} f + |\nabla f|^2| \leq \delta_k \rho_{2,k}^{-2}$. Moreover, the function $f$ satisfies $\frac{1}{2K} \, (\rho_{2,k}^{-1} \, d_{g(t)}(p_{2,t},x)+1) \leq f(x,t)+1 \leq 2K \, (\rho_{2,k}^{-1} \, d_{g(t)}(p_{2,t},x)+1)$ for all points $(x,t) \in B_{g(t_k)}(p_{2,t_k},\delta_k^{-1} \rho_{2,k}) \times [t_k-\delta_k^{-1} \rho_{2,k}^2,t_k]$.
\end{itemize}
Here, $K$ is a numerical constant. Finally, by a suitable choice of $\delta_k$, we can arrange that for each $t \in (-\infty,t_k]$ the balls $B_{g(t)}(p_{1,t},\delta_k^{-2} R(p_{1,t},t)^{-\frac{1}{2}})$ and $B_{g(t)}(p_{2,t},\delta_k^{-2} R(p_{2,t},t)^{-\frac{1}{2}})$ are disjoint.
\end{corollary}

\textbf{Proof.} 
All statements except for the last one follow immediately from Proposition \ref{cap.is.close.to.bryant}. The last statement follows from Proposition \ref{distance.between.caps}. \\

\begin{proposition}
\label{iteration}
If $(\bar{x},\bar{t}) \in M \times [t_k-2^{-j} \delta_k^{-1} \rho_{1,k}^2,t_k]$ satisfies $2^{\frac{j}{400}} \, \Lambda \rho_{1,k} \leq d_{g(\bar{t})}(p_{1,\bar{t}},\bar{x}) \leq (400KL)^{-j} \, \delta_k^{-1} \, \rho_{1,k}$, then $(\bar{x},\bar{t})$ is $2^{-j-1} \varepsilon_k$-symmetric. Similarly, if $(\bar{x},\bar{t}) \in M \times [t_k-2^{-j} \delta_k^{-1} \rho_{2,k}^2,t_k]$ satisfies $2^{\frac{j}{400}} \, \Lambda \rho_{2,k} \leq d_{g(\bar{t})}(p_{2,\bar{t}},\bar{x}) \leq (400KL)^{-j} \, \delta_k^{-1} \, \rho_{2,k}$, then $(\bar{x},\bar{t})$ is $2^{-j-1} \varepsilon_k$-symmetric. 
\end{proposition}

\textbf{Proof.} 
The proof is by induction on $j$. We first verify the assertion for $j=0$. Suppose that $(\bar{x},\bar{t}) \in M \times [t_k-2^{-j} \delta_k^{-1} \rho_{1,k}^2,t_k]$ satisfies $\Lambda \rho_{1,k} \leq d_{g(\bar{t})}(p_{1,\bar{t}},\bar{x}) \leq \delta_k^{-1} \, \rho_{1,k}$. By Corollary \ref{almost.Bryant}, the balls $B_{g(\bar{t})}(p_{1,\bar{t}},\delta_k^{-2} R(p_{1,\bar{t}},\bar{t})^{-\frac{1}{2}})$ and $B_{g(\bar{t})}(p_{2,\bar{t}},\delta_k^{-2} R(p_{2,t},\bar{t})^{-\frac{1}{2}})$ are disjoint. Since $d_{g(\bar{t})}(p_{1,\bar{t}},\bar{x}) \leq \delta_k^{-1} \, \rho_{1,k}$, we must have $d_{g(\bar{t})}(p_{2,\bar{t}},\bar{x}) \geq \delta_k^{-1} \, \rho_{2,k}$. Hence, Lemma \ref{improved.symmetry.for.faraway.points} implies that $(\bar{x},\bar{t})$ is $\frac{\varepsilon}{2}$-symmetric. This proves the assertion for $j=0$. The inductive step follows by an application of the Neck Improvement Theorem. The argument is the same as the proof of Proposition 9.11 in \cite{Brendle2}. \\

\begin{proposition}
\label{improvement.of.symmetry.of.cap}
If $k$ is sufficiently large, then $(M,g(t_k))$ is $\frac{\varepsilon_k}{2}$-symmetric. 
\end{proposition}

\textbf{Proof.} 
The arguments in Section 9 of \cite{Brendle2} go through unchanged. \\

Proposition \ref{improvement.of.symmetry.of.cap} contradicts the definition of $t_k$. This completes the proof of Theorem \ref{main.thm}.


\begin{thebibliography}{9}
\bibitem{Anderson-Chow}
G.~Anderson and B.~Chow, \textit{A pinching estimate for solutions of the linearized Ricci flow system on $3$-manifolds,} Calc. Var. PDE 23, 1--12 (2005)

\bibitem{Bamler-Kleiner}
R.~Bamler and B.~Kleiner, \textit{On the rotational symmetry of $3$-dimensional $\kappa$-solutions,} arxiv:1904.05388

\bibitem{Brendle1}
S.~Brendle, \textit{Rotational symmetry of self-similar solutions to the Ricci flow,} Invent. Math. 194, 731--764 (2013)

\bibitem{Brendle2}
S.~Brendle, \textit{Ancient solutions to the Ricci flow in dimension $3$,} arxiv:1811.02559

\bibitem{Perelman1} 
G.~Perelman, \textit{The entropy formula for the Ricci flow and its geometric applications,} arxiv:0211159

\bibitem{Perelman2}
G.~Perelman, \textit{Ricci flow with surgery on three-manifolds,} arxiv:0303109
\end{thebibliography}
\end{document}